\documentclass[11pt,reqno]{amsart}
\usepackage{graphicx}
\usepackage{verbatim}
\usepackage{textcomp}
\usepackage{amssymb}
\usepackage{cite}
\usepackage{amsmath}
\usepackage{latexsym}
\usepackage{amscd}
\usepackage{amsthm}
\usepackage{mathrsfs}
\usepackage{xypic}
\usepackage{bm}
\usepackage{url}
\usepackage{hyperref}

\newtheorem{thm}{Theorem}[section]

\newtheorem{lem}[thm]{Lemma}

\theoremstyle{definition}

\theoremstyle{remark}
\newtheorem{rem}{Remark}[section]
\numberwithin{equation}{section}
\begin{document}
\title[Some rigidity characterizations of Einstein metrics]
{Some rigidity characterizations of Einstein metrics as critical points for quadratic curvature
functionals}

\author{Bingqing Ma}
\author{Guangyue Huang}
\address{College of Mathematics and Information Science, Henan Normal
University, Xinxiang 453007, P.R. China}
\email{bqma@henannu.edu.cn (B. Ma) }
\email{hgy@henannu.edu.cn (G. Huang) }

\author{Jie Yang}
\address{College of Mathematics and System Science, Xinjiang
University, Xinjiang 830046, P.R. China} \email{yangjie1106@xju.edu.cn (J. Yang) }

\thanks{The research of author is supported by NSFC (Nos. 11371018, 11671121).}

\maketitle

\begin{abstract}

We study rigidity results for the Einstein metrics as the critical points of a family of known quadratic curvature functionals involving the scalar curvature, the Ricci curvature and the Riemannian curvature tensor, characterized by some pointwise inequalities involving the Weyl curvature and the traceless Ricci curvature. Moreover, we also provide a few rigidity results for locally conformally flat critical metrics.

\end{abstract}



\section{Introduction}

A well-known example of a Riemannian functional is
the Einstein-Hilbert functional
\begin{equation*}
\mathcal{H}=\int_{M}R
\end{equation*}
on $\mathscr{M}_1(M^{n})$, where $R$ denotes the scalar curvature and $\mathscr{M}_1(M^{n})$ is the space of equivalence classes of smooth Riemannian metrics of volume one on closed Riemannian manifold $M^{n}$, $n\geq3$. Furthermore, it is easy to see that Einstein metrics are critical for the functional $\mathcal{H}$ (see \cite{Bess2008}). In this paper, we are interested in studying the functional
\begin{equation}\label{1-Sec-1}
\mathcal{F}_{t,s}(g)=\int_M |Ric|^2+t\int_M
R^2+s\int_M |{\rm Rm}|^2
\end{equation}
where $t,s$ are real constants, $Ric$ and $Rm$ denote the Ricci curvature and the Riemannian curvature tensor, respectively. It is easy to observe from \eqref{2-Sec-2} that every Einstein metric is critical for $\mathcal{F}_{t,0}$. In \cite{Catino2015}, Catino considered the curvature functional $\mathcal{F}_{t,0}$ and obtained some conditions on the geometry of $M^n$ such that critical metrics of $\mathcal{F}_{t,0}$ are Einstein. Certainly, there exist critical metrics which are not necessarily Einstein (for
instance, see \cite[Chapter 4]{Bess2008} and \cite{Lamontagne1998}). For some
development in this direction, see
\cite{HL2004,Gursky2015,Gursky2001,Lamontagne1994,Anderson1997,Hu2010,Catino2016,HC2017,Sheng-Wang2018}
and the references therein. Therefore, it is natural to ask that under what conditions a critical metric for the functionals
$\mathcal{F}_{t,s}(s\neq0)$ must be an Einstein one.

The authors in \cite{Barros2017} show that locally conformally flat critical metrics for $\mathcal{F}_{t,s}$ with $n+4(n-1)t+4s=0$ (when $n\geq5$) and some additional conditions are space form metrics (see \cite[Theorem 3]{Barros2017}).
In this paper, we give some new characterizations, by some pointwise inequalities involving the Weyl curvature and the traceless Ricci curvature, on critical metrics for $\mathcal{F}_{t,s}$ on $\mathscr{M}_1(M^{n})$ with $n+4(n-1)t+4s\neq0$. In order to state our results, throughout this paper, we denote by $\mathring{\rm Ric}$ and $W$ the traceless Ricci tensor and the Weyl curvature, respectively.

Our main results are stated as follows:

\begin{thm}\label{4thm1}
Let $M^{n}$ be a closed manifold of dimension $n\geq3$ with
positive scalar curvature and $g$ be a critical metric for
$\mathcal{F}_{t,s}$ on $\mathscr{M}_1(M^{n})$.
Suppose that
\begin{align}\label{4-Proof-4}
\Big|W-\frac{(n-4)[4s+(n-2)]}{\sqrt{2n}(n-2)(8s+n-2)}&\mathring{\rm Ric} \mathbin{\bigcirc\mkern-15mu\wedge} g\Big||\mathring{\rm R}_{ij}|+\sqrt{\frac{2(n-1)^2}{n(n-2)}}\Big|\frac{2(n-2)s}{8s+n-2}\Big||W|^2\notag\\
<-&\sqrt{\frac{2(n-2)}{n-1}}\frac{[3n-4+2n(n-1)t+8s]}{n(8s+n-2)}R|\mathring{\rm R}_{ij}|,
\end{align}
where $t,s$ satisfy any case of the following:

(1) when $n=3,4$,
\begin{equation}\label{4th-Form-1}
\begin{cases}s>-\frac{n-2}{8}\\
n+4(n-1)t+4s\leq0\\
3n-4+2n(n-1)t+8s<0
\end{cases}
\end{equation}
or
\begin{equation}\label{4th-Form-2}
\begin{cases}
s\leq-\frac{1}{4},\,\,\quad {\rm if}\ n=3\\
s<-\frac{1}{4},\,\,\quad {\rm if}\ n=4\\
n+4(n-1)t+4s\geq0\\
3n-4+2n(n-1)t+8s>0;
\end{cases}
\end{equation}

(2) when $n\geq5$,
\begin{equation}\label{44th-Form-1}
\begin{cases}
s\geq-\frac{1}{4}\\
n+4(n-1)t+4s\leq0\\
3n-4+2n(n-1)t+8s<0
\end{cases}
\end{equation}
or
\begin{equation}\label{44th-Form-2}
\begin{cases}
s<-\frac{n-2}{8}\\
n+4(n-1)t+4s\geq0\\
3n-4+2n(n-1)t+8s>0.
\end{cases}
\end{equation}
Then $M^{n}$ is Einstein.

\end{thm}

\begin{thm}\label{3thm1}
Let $M^{n}$ be a closed manifold of dimension $n\geq3$ with
positive scalar curvature and $g$ be a critical metric for
$\mathcal{F}_{t,s}$ on $\mathscr{M}_1(M^{n})$.
Suppose that
\begin{align}\label{3th-Proof-4}
\Big|W-\frac{(n-4)[4s+(n-2)]}{\sqrt{2n}(n-2)(8s+n-2)}&\mathring{\rm Ric} \mathbin{\bigcirc\mkern-15mu\wedge} g\Big||\mathring{\rm R}_{ij}|+\sqrt{\frac{2(n-1)^2}{n(n-2)}}\Big|\frac{2(n-2)s}{8s+n-2}\Big||W|^2\notag\\
<&\sqrt{\frac{2(n-2)}{n-1}}\frac{[3n-4+2n(n-1)t+8s]}{n(8s+n-2)}R|\mathring{\rm R}_{ij}|,
\end{align}
where $t,s$ satisfy any case of the following:

(1) when $n=3$,
\begin{equation}\label{3th-Form-1}
\begin{cases}
-\frac{1}{4}\leq s<-\frac{n-4}{8}\\
n+4(n-1)t+4s\leq0\\
3n-4+2n(n-1)t+8s<0;
\end{cases}
\end{equation}

(2) when $n\geq5$,
\begin{equation}\label{33th-Form-1}
\begin{cases}
-\frac{n-4}{8}<s\leq-\frac{1}{4}\\
n+4(n-1)t+4s\geq0\\
3n-4+2n(n-1)t+8s>0.
\end{cases}
\end{equation}
Then $M^{n}$ is Einstein.

\end{thm}

\begin{thm}\label{5thm1}
Let $M^{n}$ be a closed manifold of dimension $n\geq3$ with
positive scalar curvature and $g$ be a critical metric for
$\mathcal{F}_{t,s}$ on $\mathscr{M}_1(M^{n})$, where $1+2t+2s=0$.
Suppose that
\begin{align}\label{5-Proof-1}
\Big|W+&\frac{2s(n^2-3n+4)+2(n-2)}{\sqrt{2n}(n-2)[(n-2)+2ns]}\mathring{\rm Ric}\mathbin{\bigcirc\mkern-15mu\wedge} g\Big|+\sqrt{\frac{2(n-1)^2}{n(n-2)}}\Big|\frac{(n-2)s}{(n-2)+2ns}\Big||W|^2\notag\\
\leq&\sqrt{\frac{2(n-2)}{n-1}}\frac{[2-n-n(n-1)t+2(n-2)s]}{n[(n-2)+2ns]}R|\mathring{\rm R}_{ij}|,
\end{align}
where $t,s$ satisfy any case of the following:

(1) when $n=3$,
\begin{equation}\label{5th-Form-1}
\begin{cases}
s>-\frac{1}{6}\\
2-n-n(n-1)t+2(n-2)s>0
\end{cases}
\end{equation}
or
\begin{equation}\label{5th-Form-2}
\begin{cases}
s<-\frac{1}{4}\\
2-n-n(n-1)t+2(n-2)s<0;
\end{cases}
\end{equation}

(2) when $n\geq4$,
\begin{equation}\label{55th-Form-1}
\begin{cases}
s>-\frac{1}{4}\\
2-n-n(n-1)t+2(n-2)s>0
\end{cases}
\end{equation}
or
\begin{equation}\label{55th-Form-2}
\begin{cases}
s<-\frac{n-2}{2n}\\
2-n-n(n-1)t+2(n-2)s<0.
\end{cases}
\end{equation}
Then $M^{n}$ is Einstein as long as there exists a point such that the inequality in \eqref{5-Proof-1} is strict.

\end{thm}

\begin{thm}\label{6thm1}
Let $M^{n}$ be a closed manifold of dimension $n\geq3$ with
positive scalar curvature and $g$ be a critical metric for
$\mathcal{F}_{t,s}$ on $\mathscr{M}_1(M^{n})$, where $1+2t+2s=0$.
Suppose that
\begin{align}\label{6-formProof-1}
\Big|W+&\frac{2s(n^2-3n+4)+2(n-2)}{\sqrt{2n}(n-2)[(n-2)+2ns]}\mathring{\rm Ric}\mathbin{\bigcirc\mkern-15mu\wedge} g\Big|+\sqrt{\frac{2(n-1)^2}{n(n-2)}}\Big|\frac{(n-2)s}{(n-2)+2ns}\Big||W|^2\notag\\
\leq&-\sqrt{\frac{2(n-2)}{n-1}}\frac{[2-n-n(n-1)t+2(n-2)s]}{n[(n-2)+2ns]}R|\mathring{\rm R}_{ij}|,
\end{align}
where $t,s$ satisfy any case of the following:

(1) when $n=3$,
\begin{equation}\label{6th-Form-1}
\begin{cases}
-\frac{1}{4}<s<-\frac{n-2}{2n}\\
2-n-n(n-1)t+2(n-2)s>0;
\end{cases}
\end{equation}

(2) when $n\geq5$,
\begin{equation}\label{6th-Form-2}
\begin{cases}
-\frac{n-2}{2n}<s<-\frac{1}{4}\\
2-n-n(n-1)t+2(n-2)s<0.
\end{cases}
\end{equation}
Then $M^{n}$ is Einstein as long as there exists a point such that the inequality in \eqref{6-formProof-1} is strict.

\end{thm}

\begin{thm}\label{7thm1}
Let $M^{n}$ be a closed manifold of dimension $n\geq3$ with
positive scalar curvature and $g$ be a critical metric for
$\mathcal{F}_{t,s}$ on $\mathscr{M}_1(M^{n})$.
Suppose that the divergence of Cotton tensor is zero (that is, $C_{ijk,i}=0$) and
\begin{align}\label{7-Proof-4}
\Big|W-\frac{(n-4)[4s+(n-2)]}{\sqrt{2n}(n-2)(8s+n-2)}&\mathring{\rm Ric} \mathbin{\bigcirc\mkern-15mu\wedge} g\Big||\mathring{\rm R}_{ij}|+\sqrt{\frac{2(n-1)^2}{n(n-2)}}\Big|\frac{2(n-2)s}{8s+n-2}\Big||W|^2\notag\\
<-&\sqrt{\frac{2(n-2)}{n-1}}\frac{[3n-4+2n(n-1)t+8s]}{n(8s+n-2)}R|\mathring{\rm R}_{ij}|,
\end{align}
where $t,s$ satisfy
\begin{equation}\label{7th-Form-1}
\begin{cases}s>-\frac{n-2}{8}\\
n+4(n-1)t+4s\leq0\\
3n-4+2n(n-1)t+8s<0
\end{cases}
\end{equation}
or
\begin{equation}\label{7th-Form-2}
\begin{cases}
s<-\frac{n-2}{8}\\
n+4(n-1)t+4s\geq0\\
3n-4+2n(n-1)t+8s>0.
\end{cases}
\end{equation}
Then $M^{n}$ is Einstein.

\end{thm}

\begin{thm}\label{8thm1}
Let $M^{n}$ be a closed manifold of dimension $n\geq3$ with
positive scalar curvature and $g$ be a critical metric for
$\mathcal{F}_{t,s}$ on $\mathscr{M}_1(M^{n})$.
Suppose that the divergence of Cotton tensor is zero (that is, $C_{ijk,i}=0$) and
\begin{align}\label{8th-Proof-4}
\Big|W-\frac{(n-4)[4s+(n-2)]}{\sqrt{2n}(n-2)(8s+n-2)}&\mathring{\rm Ric} \mathbin{\bigcirc\mkern-15mu\wedge} g\Big||\mathring{\rm R}_{ij}|+\sqrt{\frac{2(n-1)^2}{n(n-2)}}\Big|\frac{2(n-2)s}{8s+n-2}\Big||W|^2\notag\\
<&\sqrt{\frac{2(n-2)}{n-1}}\frac{[3n-4+2n(n-1)t+8s]}{n(8s+n-2)}R|\mathring{\rm R}_{ij}|,
\end{align}
where $t,s$ satisfy
\begin{equation}\label{8th-Form-1}
\begin{cases}
s<-\frac{n-4}{8}\\
n+4(n-1)t+4s\leq0\\
3n-4+2n(n-1)t+8s<0;
\end{cases}
\end{equation}
or
\begin{equation}\label{8th-Form-2}
\begin{cases}
s>-\frac{n-4}{8}\\
n+4(n-1)t+4s\geq0\\
3n-4+2n(n-1)t+8s>0.
\end{cases}
\end{equation}
Then $M^{n}$ is Einstein.

\end{thm}

\begin{thm}\label{11thm1}
Let $M^{n}$ be a locally conformally flat closed manifold of dimension $n\geq4$ with
positive scalar curvature and $g$ be a critical metric for
$\mathcal{F}_{t,s}$ on $\mathscr{M}_1(M^{n})$.

(1) If $n=4$ and $3t+s+1\neq0$, then $M^{4}$ is of positive constant sectional curvature;

(2) If $n\geq5$ and $t,s$ satisfy
\begin{equation}\label{11th-Form-1}
\begin{cases}
s\geq-\frac{n-2}{4}\\
(n-1)(n-2)t+2s+(n-2)<0\\
2n(n-1)t+4(n-2)s+(n^2-3n+4)\leq0
\end{cases}
\end{equation}
or
\begin{equation}\label{11th-Form-2}
\begin{cases}
s\leq-\frac{n-2}{4}\\
(n-1)(n-2)t+2s+(n-2)>0\\
2n(n-1)t+4(n-2)s+(n^2-3n+4)\geq0,
\end{cases}
\end{equation}
then $M^{n}$ is of positive constant sectional curvature.

\end{thm}

Next, we give some rigidity results for $n=3$:

\begin{thm}\label{9thm1}
Let $M^{3}$ be a closed manifold with
positive scalar curvature and $g$ be a critical metric for
$\mathcal{F}_{t,s}$ on $\mathscr{M}_1(M^{3})$, where $t,s$ satisfy
\begin{equation}\label{9th-Form-1}
\begin{cases}
s\geq-\frac{1}{4}\\
2t+2s+1>0\\
3t+s+1\geq0
\end{cases}
\end{equation}
or
\begin{equation}\label{9th-Form-2}
\begin{cases}
s\leq-\frac{1}{4}\\
2t+2s+1<0\\
3t+s+1\leq0.
\end{cases}
\end{equation}
Suppose that the divergence of Cotton tensor is zero (that is, $C_{ijk,i}=0$).
Then $M^{3}$ is of positive constant sectional curvature.

\end{thm}

\begin{thm}\label{10thm1}
Let $M^{3}$ be a closed manifold with
positive scalar curvature and $g$ be a critical metric for
$\mathcal{F}_{t,s}$ on $\mathscr{M}_1(M^{3})$, where $s=-\frac{1}{4}$ and $t\neq-\frac{1}{4}$. Then $M^{3}$ 
is of positive constant sectional curvature.

\end{thm}

\begin{rem}
In particular, when $n=3$, we have $W=0$ automatically. Hence it is seen from \eqref{2-Sec-4} that an Einstein manifold $M^3$ with positive scalar curvature must be of positive constant sectional curvature.
\end{rem}

\begin{rem}
When $s=0$, it is easy to check that our Theorem \ref{4thm1} and \ref{5thm1} reduce  Theorem 1.1 and 1.3 of \cite{MHLC2018}, respectively.
\end{rem}

\begin{rem}
For $n\geq4$, the Bach tensor is defined (see \cite{CaoChen1,HW2013}) by
\begin{equation}\label{Pre5}
B_{ij}=\frac{1}{n-3}W_{ikjl,lk}+\frac{1}{n-2}W_{ikjl}R^{kl}.
\end{equation}
By virtue of \eqref{2-Sec-8}, we have that \eqref{Pre5} can be written
as
\begin{equation}\label{Pre6}
B_{ij}=\frac{1}{n-2}(C_{kij,k}+W_{ikjl}R^{kl}).
\end{equation}
Therefore, we can define the Bach tensor on $M^3$ by
\begin{equation}\label{Pre7}
B_{ij}=C_{kij,k}.
\end{equation}
Thus, when $n=3$, $C_{ijk,i}=0$ is equivalent to $B_{ij}=0$. In \cite{Sheng-Wang2018}, Sheng and Wang
studied the case that the critical metrics are Bach-flat (that is, $B_{ij}=0$). Our Theorem \ref{9thm1} generalizes partially
the results of Sheng and Wang in \cite{Sheng-Wang2018}.
\end{rem}

\begin{rem}
When $n\geq5$, taking $s=-\frac{n-2}{4}$ in \eqref{11th-Form-1} and \eqref{11th-Form-2}, it is easy to obtain that for $t\neq-\frac{1}{2(n-1)}$, locally conformally flat closed manifold $M^n$ must be of positive constant sectional curvature. Hence, our Theorem \ref{11thm1} generalizes Theorem 4 of Barros and Da Silva \cite{Barros2017}. Moreover, for $n=3,4$, our Theorems \ref{11thm1} and \ref{10thm1} can be seen as a supplementary to Theorem 4 of Barros and Da Silva in \cite{Barros2017}.

\end{rem}

\section{Preliminaries}
For $n\geq3$, it is well-known that the Weyl curvature tensor and the Cotton tensor are defined by
\begin{align}\label{2-Sec-4}
W_{ijkl}=&R_{ijkl}-\frac{1}{n-2}(R_{ik}g_{jl}-R_{il}g_{jk}
+R_{jl}g_{ik}-R_{jk}g_{il})\notag\\
&+\frac{R}{(n-1)(n-2)}(g_{ik}g_{jl}-g_{il}g_{jk})\notag\\
=&R_{ijkl}-\frac{1}{n-2}(\mathring{\rm R}_{ik}g_{jl}-\mathring{\rm R}_{il}g_{jk}
+\mathring{\rm R}_{jl}g_{ik}-\mathring{\rm R}_{jk}g_{il})\notag\\
&-\frac{R}{n(n-1)}(g_{ik}g_{jl}-g_{il}g_{jk}),
\end{align}
and
\begin{align}\label{2-Sec-5}
C_{ijk}=&R_{kj,i}-R_{ki,j}-\frac{1}{2(n-1)}(R_{,i}g_{jk}-R_{,j}g_{ik})\notag\\
=&\mathring{\rm R}_{kj,i}-\mathring{\rm R}_{ki,j}+\frac{n-2}{2n(n-1)}(R_{,i}g_{jk}-R_{,j}g_{ik}),
\end{align}
respectively. Here $\mathring{\rm R}_{ij}=R_{ij}-\frac{1}{n} Rg_{ij}$ denotes the traceless Ricci tensor and the indices after a comma denote the covariant derivatives. From the definition of the Cotton tensor, it is easy to see
$$
C_{ijk}=-C_{jik},\qquad g^{ij}C_{ijk}=g^{ik}C_{ijk}=g^{jk}C_{ijk}=0
$$
and
$$
C_{ijk,k}=0,\qquad C_{ijk}+C_{jki}+C_{kij}=0.
$$
For $n\geq4$, the divergence of the Weyl curvature tensor is related to the Cotton
tensor by
\begin{equation}\label{2-Sec-8}
-\frac{n-3}{n-2}C_{ijk}=W_{ijkl,l}.
\end{equation}
Moreover, $W_{ijkl}=0$ holds naturally on $(M^3,g)$, and $(M^3,g)$ is
locally conformally flat if and only if $C_{ijk}=0$. For $n\geq4$,
$(M^n,g)$ is locally conformally flat if and only if $W_{ijkl}=0$.

It has been proved by Catino in \cite{Catino2015} (see \cite[Proposition 6.1]{Catino2015}) that
a metric $g$ is critical for
$\mathcal{F}_{t,s}$ on $\mathscr{M}_1(M^{n})$ if and only if it satisfies the equations
\begin{align}\label{2-Sec-2}
(1+4s)\Delta \mathring{\rm R}_{ij}=&(1+2t+2s)R_{,ij}-\frac{1+2t+2s}{n}(\Delta R)g_{ij}-2(1+2s)R_{ikjl}\mathring{\rm R}_{kl}\notag\\
&-\frac{2+2nt-4s}{n}R\mathring{\rm R}_{ij}+\frac{2}{n}(|\mathring{\rm R}_{ij}|^2+s|{\rm Rm}|^2)g_{ij}\notag\\
&-2sR_{ikpq}R_{jkpq}+4s\mathring{\rm R}_{ik}\mathring{\rm R}_{jk}
\end{align}
and
\begin{align}\label{2-Sec-3}
{[n+4(n-1)t+4s]}\Delta R=&(n-4)(|R_{ij}|^2+tR^2+s|{\rm Rm}|^2-\lambda)\notag\\
=&(n-4)\Big[s|W|^2+\frac{n-2+4s}{n-2}|\mathring{\rm R}_{ij}|^2\notag\\
&+\frac{n-1+n(n-1)t+2s}{n(n-1)}R^2-\lambda\Big],
\end{align}
where $\lambda=\mathcal{F}_{t,s}(g)$ and we used the fact
\begin{align}\label{2-Sec-17}
|{\rm Rm}|^2=&|W|^2+\frac{4}{n-2}|\mathring{\rm R}_{ij}|^2+\frac{2}{n(n-1)}R^2
\end{align}
from \eqref{2-Sec-4}.

Using the formula \eqref{2-Sec-4}, we can also derive
\begin{align}\label{2-Sec-15}
\mathring{\rm R}_{kl}R_{ikjl}=&\mathring{\rm R}_{kl}W_{ikjl}+\frac{1}{n-2}(|\mathring{\rm R}_{ij}|^2g_{ij}-2\mathring{\rm R}_{ik}\mathring{\rm R}_{jk})-\frac{1}{n(n-1)}R \mathring{\rm R}_{ij}
\end{align}
and
\begin{align}\label{2-Sec-16}
R_{ikpq}R_{jkpq}=&W_{ikpq}W_{jkpq}+\frac{4}{n-2}W_{ikjl}\mathring{\rm R}_{kl}+\frac{2(n-4)}{(n-2)^2}\mathring{\rm R}_{ik}\mathring{\rm R}_{jk}\notag\\
&+\frac{2}{(n-2)^2}|\mathring{\rm R}_{ij}|^2g_{ij}+\frac{2}{n^2(n-1)}R^2g_{ij}+\frac{4}{n(n-1)}R\mathring{\rm R}_{ij}.
\end{align}
Therefore, \eqref{2-Sec-2} can be written as
\begin{align}\label{add2-Sec-2}
(1+4s)\Delta \mathring{\rm R}_{ij}=&(1+2t+2s)\mathring{R}_{,ij}-2(1+2s)R_{ikjl}\mathring{\rm R}_{kl}\notag\\
&-\frac{2+2nt-4s}{n}R\mathring{\rm R}_{ij}+\frac{2}{n}(|\mathring{\rm R}_{ij}|^2+s|{\rm Rm}|^2)g_{ij}\notag\\
&-2sR_{ikpq}R_{jkpq}+4s\mathring{\rm R}_{ik}\mathring{\rm R}_{jk}\notag\\
=&(1+2t+2s)\mathring{R}_{,ij}-\frac{2(n-2)+4ns}{n-2}W_{ikjl}\mathring{\rm R}_{kl}-2sW_{ikpq}W_{jkpq}\notag\\
&+\Big[-\frac{4s(n^2-3n+4)+4(n-2)}{n(n-2)^2}|\mathring{\rm R}_{ij}|^2+\frac{2s}{n}|W|^2\Big]g_{ij}\notag\\
&+\frac{4s(n^2-3n+4)+4(n-2)}{(n-2)^2}\mathring{\rm R}_{ik}\mathring{\rm R}_{jk}\notag\\
&+\frac{4-2n-2n(n-1)t+4(n-2)s}{n(n-1)}R\mathring{\rm R}_{ij},
\end{align}
where $\mathring{R}_{,ij}=R_{,ij}-\frac{1}{n}(\Delta R)g_{ij}$. It follows from \eqref{add2-Sec-2} that
\begin{align}\label{add2-Sec-3}
\frac{1+4s}{2}\Delta |\mathring{\rm R}_{ij}|^2=&(1+4s)|\nabla \mathring{\rm R}_{ij}|^2+(1+4s)\mathring{\rm R}_{ij}\Delta \mathring{\rm R}_{ij}\notag\\
=&(1+4s)|\nabla \mathring{\rm R}_{ij}|^2+(1+2t+2s)R_{,ij}\mathring{\rm R}_{ij}\notag\\
&-\frac{2(n-2)+4ns}{n-2}W_{ikjl}\mathring{\rm R}_{kl}\mathring{\rm R}_{ij}-2sW_{ikpq}W_{jkpq}\mathring{\rm R}_{ij}\notag\\
&+\frac{4s(n^2-3n+4)+4(n-2)}{(n-2)^2}\mathring{\rm R}_{ik}\mathring{\rm R}_{kj}\mathring{\rm R}_{ji}\notag\\
&+\frac{4-2n-2n(n-1)t+4(n-2)s}{n(n-1)}R|\mathring{\rm R}_{ij}|^2.
\end{align}
Integrating both sides of \eqref{add2-Sec-3} yields
\begin{align}\label{Aug2-Proof-1}
0=&(1+4s)\int_M|\nabla \mathring{\rm R}_{ij}|^2+\int_M\Big(-\frac{(n-2)(1+2t+2s)}{2n}|\nabla R|^2\notag\\
&-\frac{2(n-2)+4ns}{n-2}W_{ikjl}\mathring{\rm R}_{kl}\mathring{\rm R}_{ij}-2sW_{ikpq}W_{jkpq}\mathring{\rm R}_{ij}\notag\\
&+\frac{4s(n^2-3n+4)+4(n-2)}{(n-2)^2}\mathring{\rm R}_{ik}\mathring{\rm R}_{kj}\mathring{\rm R}_{ji}\notag\\
&+\frac{4-2n-2n(n-1)t+4(n-2)s}{n(n-1)}R|\mathring{\rm R}_{ij}|^2\Big),
\end{align}
where we used the second Bianchi identity $\mathring{\rm R}_{kj,k}=\frac{n-2}{2n}R_{,j}$.
Hence, we obtain the following result:

\begin{lem}\label{4-Lemma4}
Let $M^n$ be a closed manifold and $g$ be a critical metric for
$\mathcal{F}_{t,s}$ on $\mathscr{M}_1(M^{n})$. Then
\begin{align}\label{Aug2-Proof-2}
(1+4s)\int_M|\nabla \mathring{\rm R}_{ij}|^2=&\int_M\Big(\frac{(n-2)(1+2t+2s)}{2n}|\nabla R|^2\notag\\
&+\frac{2(n-2)+4ns}{n-2}W_{ikjl}\mathring{\rm R}_{kl}\mathring{\rm R}_{ij}+2sW_{ikpq}W_{jkpq}\mathring{\rm R}_{ij}\notag\\
&-\frac{4s(n^2-3n+4)+4(n-2)}{(n-2)^2}\mathring{\rm R}_{ik}\mathring{\rm R}_{kj}\mathring{\rm R}_{ji}\notag\\
&-\frac{4-2n-2n(n-1)t+4(n-2)s}{n(n-1)}R|\mathring{\rm R}_{ij}|^2\Big).
\end{align}

\end{lem}

For any closed manifold, we also have the following result (see \cite[Lemma 2.2]{MHLC2018})

\begin{lem}\label{4-Lemma1}
Let $M^n$ be a closed manifold. Then
\begin{align}\label{4-Proof-1}
\int_M|\nabla \mathring{\rm R}_{ij}|^2=&\int_M\Big(
W_{ijkl}\mathring{\rm R}_{jl}\mathring{\rm R}_{ik}
-\frac{n}{n-2}\mathring{\rm R}_{ij}\mathring{\rm R}_{jk}
\mathring{\rm R}_{ki}\notag\\
&-\frac{1}{n-1}R|\mathring{\rm R}_{ij}|^2+\frac{(n-2)^2}{4n(n-1)}|\nabla R|^2+\frac{1}{2}|C_{ijk}|^2\Big).
\end{align}

\end{lem}

The next lemma comes from \cite{HM-2016,FuPeng-2017Hokkaido,MH2018} (for the case of $\lambda=\frac{2}{n-2}$, see \cite{Catino-adv2016}):

\begin{lem}\label{4-Lemma2}
For every Riemannian manifold $(M^n,g)$ and any $\lambda\in \mathbb{R}$, the following estimate holds
\begin{align}\label{4-Proof-2}
\Big|-W_{ijkl}\mathring{\rm R}_{jl}&\mathring{\rm R}_{ik}
+\lambda \mathring{\rm R}_{ij}\mathring{\rm R}_{jk}\mathring{\rm R}_{ki}\Big|\notag\\
\leq&
\sqrt{\frac{n-2}{2(n-1)}}\Big(|W|^2+\frac{2(n-2)\lambda^2}{n}|\mathring{\rm R}_{ij}|^2
\Big)^{\frac{1}{2}}|\mathring{\rm R}_{ij}|^2\notag\\
=&\sqrt{\frac{n-2}{2(n-1)}}\Big|W+\frac{\lambda}{\sqrt{2n}}\mathring{\rm Ric} \mathbin{\bigcirc\mkern-15mu\wedge} g\Big||\mathring{\rm R}_{ij}|^2.
\end{align}

\end{lem}

\section{Proof of main results}

\subsection{Proof of Theorem \ref{4thm1}}
Notice that \eqref{4-Proof-1} can be written as
\begin{align}\label{4-Proof-7}
(1+4s)\int_M|\nabla \mathring{\rm R}_{ij}|^2=&(1+4s)\int_M\Big(
W_{ijkl}\mathring{\rm R}_{jl}\mathring{\rm R}_{ik}
-\frac{n}{n-2}\mathring{\rm R}_{ij}\mathring{\rm R}_{jk}
\mathring{\rm R}_{ki}\notag\\
&-\frac{1}{n-1}R|\mathring{\rm R}_{ij}|^2+\frac{(n-2)^2}{4n(n-1)}|\nabla R|^2+\frac{1}{2}|C_{ijk}|^2\Big).
\end{align}
Combining \eqref{4-Proof-7} with \eqref{Aug2-Proof-2}, we have
\begin{align}\label{4-Proof-8}
0=&\int_M\Big[ \frac{n-2+8s}{n-2}
W_{ijkl}\mathring{\rm R}_{jl}\mathring{\rm R}_{ik}
+\frac{(n-4)[4s+(n-2)]}{(n-2)^2}\mathring{\rm R}_{ij}\mathring{\rm R}_{jk}
\mathring{\rm R}_{ki}\notag\\
&+2sW_{ikpq}W_{jkpq}\mathring{\rm R}_{ij}+\frac{3n-4+2n(n-1)t+8s}{n(n-1)}R|\mathring{\rm R}_{ij}|^2\notag\\
&+\frac{(n-2)[n+4(n-1)t+4s]}{4n(n-1)}|\nabla R|^2-\frac{1+4s}{2}|C_{ijk}|^2\Big],
\end{align}
which is equivalent to
\begin{align}\label{4-Proof-9}
0=&\int_M\Big[-W_{ijkl}\mathring{\rm R}_{jl}\mathring{\rm R}_{ik}
-\frac{(n-4)[4s+(n-2)]}{(n-2)(8s+n-2)}\mathring{\rm R}_{ij}\mathring{\rm R}_{jk}
\mathring{\rm R}_{ki}\notag\\
&-\frac{2(n-2)s}{8s+n-2}W_{ikpq}W_{jkpq}\mathring{\rm R}_{ij}-\frac{(n-2)[3n-4+2n(n-1)t+8s]}{n(n-1)(8s+n-2)}R|\mathring{\rm R}_{ij}|^2\notag\\
&-\frac{(n-2)^2[n+4(n-1)t+4s]}{4n(n-1)(8s+n-2)}|\nabla R|^2+\frac{(n-2)(1+4s)}{2(8s+n-2)}|C_{ijk}|^2\Big]
\end{align}
as long as $8s+n-2\neq0$.
Substituting the estimate \eqref{4-Proof-2} with $\lambda=-\frac{(n-4)[4s+(n-2)]}{(n-2)(8s+n-2)}$ and
\begin{align}\label{add4-Proof-10}
|W_{ikpq}W_{jkpq}\mathring{\rm R}_{ij}|\leq\sqrt{\frac{n-1}{n}}|W|^2|\mathring{\rm R}_{ij}|
\end{align}
into \eqref{4-Proof-9} gives
\begin{align}\label{4-Proof-10}
0\geq&\int_M\Bigg[-\sqrt{\frac{n-2}{2(n-1)}}\Big|W-\frac{(n-4)[4s+(n-2)]}{\sqrt{2n}(n-2)(8s+n-2)}\mathring{\rm Ric}  \mathbin{\bigcirc\mkern-15mu\wedge} g\Big||\mathring{\rm R}_{ij}|\notag\\
&-\sqrt{\frac{n-1}{n}}\Big|\frac{2(n-2)s}{8s+n-2}\Big||W|^2-\frac{(n-2)[3n-4+2n(n-1)t+8s]}{n(n-1)(8s+n-2)}R|\mathring{\rm R}_{ij}|\Bigg]|\mathring{\rm R}_{ij}|\notag\\
&+\int_M\Bigg[-\frac{(n-2)^2[n+4(n-1)t+4s]}{4n(n-1)(8s+n-2)}|\nabla R|^2+\frac{(n-2)(1+4s)}{2(8s+n-2)}|C_{ijk}|^2\Bigg].
\end{align}
For the proof of \eqref{add4-Proof-10}, we refer to \cite[Lemma 2.4]{Huisken1985}. Noticing that if $t,s$ satisfy \eqref{4th-Form-1} or \eqref{44th-Form-1}, then we have
\begin{equation}\label{4th-Form-13}
\begin{cases}
1+4s\geq0\\
8s+n-2>0\\
n+4(n-1)t+4s\leq0\\
3n-4+2n(n-1)t+8s<0.
\end{cases}
\end{equation}
Therefore, applying \eqref{4-Proof-4} and \eqref{4th-Form-13} into \eqref{4-Proof-10} gives
\begin{align}\label{4-Proof-11}
0\geq&\int_M\Bigg[-\sqrt{\frac{n-2}{2(n-1)}}\Big|W-\frac{(n-4)[4s+(n-2)]}{\sqrt{2n}(n-2)(8s+n-2)}\mathring{\rm Ric} \mathbin{\bigcirc\mkern-15mu\wedge} g\Big||\mathring{\rm R}_{ij}|\notag\\
&-\sqrt{\frac{n-1}{n}}\Big|\frac{2(n-2)s}{8s+n-2}\Big||W|^2-\frac{(n-2)[3n-4+2n(n-1)t+8s]}{n(n-1)(8s+n-2)}R|\mathring{\rm R}_{ij}|\Bigg]|\mathring{\rm R}_{ij}|\notag\\
&+\int_M\Bigg[-\frac{(n-2)^2[n+4(n-1)t+4s]}{4n(n-1)(8s+n-2)}|\nabla R|^2+\frac{(n-2)(1+4s)}{2(8s+n-2)}|C_{ijk}|^2\Bigg]\notag\\
\geq&0,
\end{align}
which shows $\mathring{\rm R}_{ij}=0$ and hence $M^n$ is Einstein.

Similarly, if $t,s$ satisfy \eqref{4th-Form-2} or \eqref{44th-Form-2}, then we have
\begin{equation}\label{4th-Form-14}
\begin{cases}
1+4s\leq0\\
8s+n-2<0\\
n+4(n-1)t+4s\geq0\\
3n-4+2n(n-1)t+8s>0.
\end{cases}
\end{equation}
Therefore, applying \eqref{4-Proof-4} and \eqref{4th-Form-14} into \eqref{4-Proof-10} also yields the estimate \eqref{4-Proof-11}
 and the desired Theorem \ref{4thm1} follows.

\subsection{Proof of Theorem \ref{3thm1}}
When $8s+n-2\neq0$, inserting the estimate \eqref{4-Proof-2} with $\lambda=-\frac{(n-4)[4s+(n-2)]}{(n-2)(8s+n-2)}$ and
\eqref{add4-Proof-10} into \eqref{4-Proof-9}, we deduce
\begin{align}\label{3-Proof-1}
0\leq&\int_M\Bigg[\sqrt{\frac{n-2}{2(n-1)}}\Big|W-\frac{(n-4)[4s+(n-2)]}{\sqrt{2n}(n-2)(8s+n-2)}\mathring{\rm Ric}  \mathbin{\bigcirc\mkern-15mu\wedge} g\Big||\mathring{\rm R}_{ij}|\notag\\
&+\sqrt{\frac{n-1}{n}}\Big|\frac{2(n-2)s}{8s+n-2}\Big||W|^2-\frac{(n-2)[3n-4+2n(n-1)t+8s]}{n(n-1)(8s+n-2)}R|\mathring{\rm R}_{ij}|\Bigg]|\mathring{\rm R}_{ij}|\notag\\
&+\int_M\Bigg[-\frac{(n-2)^2[n+4(n-1)t+4s]}{4n(n-1)(8s+n-2)}|\nabla R|^2+\frac{(n-2)(1+4s)}{2(8s+n-2)}|C_{ijk}|^2\Bigg].
\end{align}
When $n=3$, if $t,s$ satisfy \eqref{3th-Form-1}, then we have
\begin{equation}\label{3th-Proof-2}
\begin{cases}
1+4s\geq0\\
8s+n-2<0\\
n+4(n-1)t+4s\leq0\\
3n-4+2n(n-1)t+8s<0.
\end{cases}
\end{equation}
Therefore, applying \eqref{3th-Proof-2} and \eqref{3th-Proof-4} into \eqref{3-Proof-1} gives
\begin{align}\label{3-Proof-3}
0\leq&\int_M\Bigg[\sqrt{\frac{n-2}{2(n-1)}}\Big|W-\frac{(n-4)[4s+(n-2)]}{\sqrt{2n}(n-2)(8s+n-2)}\mathring{\rm Ric}  \mathbin{\bigcirc\mkern-15mu\wedge} g\Big||\mathring{\rm R}_{ij}|\notag\\
&+\sqrt{\frac{n-1}{n}}\Big|\frac{2(n-2)s}{8s+n-2}\Big||W|^2-\frac{(n-2)[3n-4+2n(n-1)t+8s]}{n(n-1)(8s+n-2)}R|\mathring{\rm R}_{ij}|\Bigg]|\mathring{\rm R}_{ij}|\notag\\
&+\int_M\Bigg[-\frac{(n-2)^2[n+4(n-1)t+4s]}{4n(n-1)(8s+n-2)}|\nabla R|^2+\frac{(n-2)(1+4s)}{2(8s+n-2)}|C_{ijk}|^2\Bigg]\notag\\
\leq&0,
\end{align}
which shows that $M^3$ is Einstein.

On the other hand, if $t,s$ satisfy \eqref{33th-Form-1}, then we have
\begin{equation}\label{3-Proof-4}
\begin{cases}
1+4s\leq0\\
8s+n-2>0\\
n+4(n-1)t+4s\geq0\\
3n-4+2n(n-1)t+8s>0.
\end{cases}
\end{equation}
Applying \eqref{3-Proof-4} and \eqref{3th-Proof-4} into \eqref{3-Proof-1} also yields the same estimate \eqref{3-Proof-3}
and the desired Theorem \ref{3thm1} follows.

\subsection{Proof of Theorem \ref{5thm1}}
When $t,s$ satisfy $1+2t+2s=0$, then the formula \eqref{add2-Sec-3} becomes
\begin{align}\label{5-Proof-2}
\frac{1+4s}{2}\Delta |\mathring{\rm R}_{ij}|^2=&(1+4s)|\nabla \mathring{\rm R}_{ij}|^2-\frac{2(n-2)+4ns}{n-2}W_{ikjl}\mathring{\rm R}_{kl}\mathring{\rm R}_{ij}\notag\\
&-2sW_{ikpq}W_{jkpq}\mathring{\rm R}_{ij}+\frac{4s(n^2-3n+4)+4(n-2)}{(n-2)^2}\mathring{\rm R}_{ik}\mathring{\rm R}_{kj}\mathring{\rm R}_{ji}\notag\\
&+\frac{4-2n-2n(n-1)t+4(n-2)s}{n(n-1)}R|\mathring{\rm R}_{ij}|^2,
\end{align}
which gives
\begin{align}\label{5-Proof-3}
&\frac{(n-2)(1+4s)}{4[(n-2)+2ns]}\Delta |\mathring{\rm R}_{ij}|^2=\frac{(n-2)(1+4s)}{2[(n-2)+2ns]}|\nabla \mathring{\rm R}_{ij}|^2-W_{ikjl}\mathring{\rm R}_{kl}\mathring{\rm R}_{ij}\notag\\
&-\frac{(n-2)s}{(n-2)+2ns}W_{ikpq}W_{jkpq}\mathring{\rm R}_{ij}+\frac{2s(n^2-3n+4)+2(n-2)}{(n-2)[(n-2)+2ns]}\mathring{\rm R}_{ik}\mathring{\rm R}_{kj}\mathring{\rm R}_{ji}\notag\\
&+\frac{(n-2)[2-n-n(n-1)t+2(n-2)s]}{n(n-1)[(n-2)+2ns]}R|\mathring{\rm R}_{ij}|^2\notag\\
\geq&\frac{(n-2)(1+4s)}{2[(n-2)+2ns]}|\nabla \mathring{\rm R}_{ij}|^2-\sqrt{\frac{n-2}{2(n-1)}}\Big|W+\frac{2s(n^2-3n+4)+2(n-2)}{\sqrt{2n}(n-2)[(n-2)+2ns]}\mathring{\rm Ric} \mathbin{\bigcirc\mkern-15mu\wedge} g\Big||\mathring{\rm R}_{ij}|^2\notag\\
&-\sqrt{\frac{n-1}{n}}\Big|\frac{(n-2)s}{(n-2)+2ns}\Big||W|^2|\mathring{\rm R}_{ij}|+\frac{(n-2)[2-n-n(n-1)t+2(n-2)s]}{n(n-1)[(n-2)+2ns]}R|\mathring{\rm R}_{ij}|^2
\end{align}
provided $n-2+2ns\neq0$. Noticing that if $t,s$ satisfy \eqref{5th-Form-1} or \eqref{55th-Form-1}, then we have
\begin{equation}\label{5-Proof-4}
\begin{cases}
1+4s>0\\
n-2+2ns>0\\
2-n-n(n-1)t+2n(n-2)s>0.
\end{cases}
\end{equation}
Similarly, if $t,s$ satisfy \eqref{5th-Form-2} or \eqref{55th-Form-2}, then we have
\begin{equation}\label{5-Proof-5}
\begin{cases}
1+4s<0\\
n-2+2ns<0\\
2-n-n(n-1)t+2n(n-2)s<0.
\end{cases}
\end{equation}
Clearly, if \eqref{5-Proof-4} or \eqref{5-Proof-5} holds, then from \eqref{5-Proof-3} and \eqref{5-Proof-1} we both have
\begin{align}\label{5-Proof-6}
&\frac{(n-2)(1+4s)}{4[(n-2)+2ns]}\Delta |\mathring{\rm R}_{ij}|^2\notag\\
\geq&\frac{(n-2)(1+4s)}{2[(n-2)+2ns]}|\nabla \mathring{\rm R}_{ij}|^2-\sqrt{\frac{n-2}{2(n-1)}}\Big|W+\frac{2s(n^2-3n+4)+2(n-2)}{\sqrt{2n}(n-2)[(n-2)+2ns]}\mathring{\rm Ric} \mathbin{\bigcirc\mkern-15mu\wedge} g\Big||\mathring{\rm R}_{ij}|^2\notag\\
&-\sqrt{\frac{n-1}{n}}\Big|\frac{(n-2)s}{(n-2)+2ns}\Big||W|^2|\mathring{\rm R}_{ij}|+\frac{(n-2)[2-n-n(n-1)t+2(n-2)s]}{n(n-1)[(n-2)+2ns]}R|\mathring{\rm R}_{ij}|^2\notag\\
\geq&0,
\end{align}
which shows that $|\mathring{\rm R}_{ij}|^2$ is subharmonic on $M^n$. Using the maximum principle, we obtain that $|\mathring{\rm R}_{ij}|$ is constant and $\nabla \mathring{\rm R}_{ij}=0$, implying that the Ricci curvature is parallel and the scalar curvature $R$ is constant. In particular, \eqref{5-Proof-3} becomes
\begin{align}\label{5-Proof-7}
0=&-\sqrt{\frac{n-2}{2(n-1)}}\Big|W+\frac{2s(n^2-3n+4)+2(n-2)}{\sqrt{2n}(n-2)[(n-2)+2ns]}\mathring{\rm Ric} \mathbin{\bigcirc\mkern-15mu\wedge} g\Big||\mathring{\rm R}_{ij}|\notag\\
&-\sqrt{\frac{n-1}{n}}\Big|\frac{(n-2)s}{(n-2)+2ns}\Big||W|^2+\frac{(n-2)[2-n-n(n-1)t+2(n-2)s]}{n(n-1)[(n-2)+2ns]}R|\mathring{\rm R}_{ij}|.
\end{align}
If there exists a point $p$ such that the inequality \eqref{5-Proof-1} is strict, then from \eqref{5-Proof-7} we have
$|\mathring{\rm R}_{ij}|(p)=0$ which with the fact that $|\mathring{\rm R}_{ij}|$ constant shows that $\mathring{\rm R}_{ij}=0$, that is, $M^n$ is Einstein, completing the proof of Theorem \ref{5thm1}.

\subsection{Proof of Theorem \ref{6thm1}}
When $1+2t+2s=0$, \eqref{5-Proof-2} can also be written as
\begin{align}\label{6-Proof-1}
&-\frac{(n-2)(1+4s)}{4[(n-2)+2ns]}\Delta |\mathring{\rm R}_{ij}|^2=-\frac{(n-2)(1+4s)}{2[(n-2)+2ns]}|\nabla \mathring{\rm R}_{ij}|^2+W_{ikjl}\mathring{\rm R}_{kl}\mathring{\rm R}_{ij}\notag\\
&+\frac{(n-2)s}{(n-2)+2ns}W_{ikpq}W_{jkpq}\mathring{\rm R}_{ij}-\frac{2s(n^2-3n+4)+2(n-2)}{(n-2)[(n-2)+2ns]}\mathring{\rm R}_{ik}\mathring{\rm R}_{kj}\mathring{\rm R}_{ji}\notag\\
&-\frac{(n-2)[2-n-n(n-1)t+2(n-2)s]}{n(n-1)[(n-2)+2ns]}R|\mathring{\rm R}_{ij}|^2.
\end{align}
Thus, we obtain
\begin{align}\label{6-Proof-2}
&-\frac{(n-2)(1+4s)}{4[(n-2)+2ns]}\Delta |\mathring{\rm R}_{ij}|^2\notag\\
\geq&-\frac{(n-2)(1+4s)}{2[(n-2)+2ns]}|\nabla \mathring{\rm R}_{ij}|^2-\sqrt{\frac{n-2}{2(n-1)}}\Big|W+\frac{2s(n^2-3n+4)+2(n-2)}{\sqrt{2n}(n-2)[(n-2)+2ns]}\mathring{\rm Ric} \mathbin{\bigcirc\mkern-15mu\wedge} g\Big||\mathring{\rm R}_{ij}|^2\notag\\
&-\sqrt{\frac{n-1}{n}}\Big|\frac{(n-2)s}{(n-2)+2ns}\Big||W|^2|\mathring{\rm R}_{ij}|-\frac{(n-2)[2-n-n(n-1)t+2(n-2)s]}{n(n-1)[(n-2)+2ns]}R|\mathring{\rm R}_{ij}|^2.
\end{align}
When $n=3$, if $t,s$ satisfy \eqref{6th-Form-1}, then we have
\begin{equation}\label{6-Proof-3}
\begin{cases}
1+4s>0\\
n-2+2ns<0\\
2-n-n(n-1)t+2(n-2)s>0.
\end{cases}
\end{equation}
Therefore, applying \eqref{3-Proof-4} and \eqref{6-formProof-1} into \eqref{6-Proof-2} yields
\begin{align}\label{6-Proof-4}
&-\frac{(n-2)(1+4s)}{4[(n-2)+2ns]}\Delta |\mathring{\rm R}_{ij}|^2\notag\\
\geq&-\frac{(n-2)(1+4s)}{2[(n-2)+2ns]}|\nabla \mathring{\rm R}_{ij}|^2-\sqrt{\frac{n-2}{2(n-1)}}\Big|W+\frac{2s(n^2-3n+4)+2(n-2)}{\sqrt{2n}(n-2)[(n-2)+2ns]}\mathring{\rm Ric} \mathbin{\bigcirc\mkern-15mu\wedge} g\Big||\mathring{\rm R}_{ij}|^2\notag\\
&-\sqrt{\frac{n-1}{n}}\Big|\frac{(n-2)s}{(n-2)+2ns}\Big||W|^2|\mathring{\rm R}_{ij}|-\frac{(n-2)[2-n-n(n-1)t+2(n-2)s]}{n(n-1)[(n-2)+2ns]}R|\mathring{\rm R}_{ij}|^2\notag\\
\geq&0,
\end{align}
which shows that $|\mathring{\rm R}_{ij}|^2$ is subharmonic on $M^3$.

When $n\geq5$, if $t,s$ satisfy \eqref{6th-Form-2}, then we have
\begin{equation}\label{6-Proof-5}
\begin{cases}
1+4s<0\\
n-2+2ns>0\\
2-n-n(n-1)t+2(n-2)s<0.
\end{cases}
\end{equation}
Therefore, applying \eqref{6-Proof-5} and \eqref{6-formProof-1} into \eqref{6-Proof-2} also yields the estimate \eqref{6-Proof-4}. Then following the proof of Theorem \ref{5thm1} line by line we finish the proof of Theorem \ref{6thm1}.

\subsection{Proof of Theorems \ref{7thm1} and \ref{8thm1}}
By the definition of the Cotton tensor given by \eqref{2-Sec-5}, we have
\begin{align}\label{7-Proof-1}
\int_MC_{ijk,i}R_{jk}=-\int_MC_{ijk}R_{jk,i}=-\frac{1}{2} \int_M|C_{ijk}|^2,
\end{align}
which shows that if $C_{ijk,i}=0$, then we have $C_{ijk}=0$. Thus, \eqref{4-Proof-10} becomes
\begin{align}\label{7-Proof-2}
0\geq&\int_M\Bigg[-\sqrt{\frac{n-2}{2(n-1)}}\Big|W-\frac{(n-4)[4s+(n-2)]}{\sqrt{2n}(n-2)(8s+n-2)}\mathring{\rm Ric}  \mathbin{\bigcirc\mkern-15mu\wedge} g\Big||\mathring{\rm R}_{ij}|\notag\\
&-\sqrt{\frac{n-1}{n}}\Big|\frac{2(n-2)s}{8s+n-2}\Big||W|^2-\frac{(n-2)[3n-4+2n(n-1)t+8s]}{n(n-1)(8s+n-2)}R|\mathring{\rm R}_{ij}|\Bigg]|\mathring{\rm R}_{ij}|\notag\\
&-\frac{(n-2)^2[n+4(n-1)t+4s]}{4n(n-1)(8s+n-2)}\int_M|\nabla R|^2.
\end{align}
If $t,s$ satisfy \eqref{7th-Form-1} or \eqref{7th-Form-2}, then $M^n$ must be Einstein and the proof of Theorem \ref{7thm1} is finished.

Similarly, if $C_{ijk,i}=0$, \eqref{3-Proof-1} becomes
\begin{align}\label{8-Proof-3}
0\leq&\int_M\Bigg[\sqrt{\frac{n-2}{2(n-1)}}\Big|W-\frac{(n-4)[4s+(n-2)]}{\sqrt{2n}(n-2)(8s+n-2)}\mathring{\rm Ric}  \mathbin{\bigcirc\mkern-15mu\wedge} g\Big||\mathring{\rm R}_{ij}|\notag\\
&+\sqrt{\frac{n-1}{n}}\Big|\frac{2(n-2)s}{8s+n-2}\Big||W|^2-\frac{(n-2)[3n-4+2n(n-1)t+8s]}{n(n-1)(8s+n-2)}R|\mathring{\rm R}_{ij}|\Bigg]|\mathring{\rm R}_{ij}|\notag\\
&-\frac{(n-2)^2[n+4(n-1)t+4s]}{4n(n-1)(8s+n-2)}\int_M|\nabla R|^2,
\end{align}
which shows that $M^n$ is Einstein as long as $t,s$ satisfy \eqref{8th-Form-1} or \eqref{8th-Form-2}. It concludes the proof of
Theorem \ref{8thm1}.

\subsection{Proof of Theorem \ref{11thm1}}
When $W_{ijkl}=0$, \eqref{Aug2-Proof-2} becomes
\begin{align}\label{10-Proof-1}
(1+4s)\int_M|\nabla \mathring{\rm R}_{ij}|^2=&\int_M\Big(\frac{(n-2)(1+2t+2s)}{2n}|\nabla R|^2\notag\\
&-\frac{4s(n^2-3n+4)+4(n-2)}{(n-2)^2}\mathring{\rm R}_{ik}\mathring{\rm R}_{kj}\mathring{\rm R}_{ji}\notag\\
&-\frac{4-2n-2n(n-1)t+4(n-2)s}{n(n-1)}R|\mathring{\rm R}_{ij}|^2\Big)
\end{align}
and \eqref{4-Proof-1} becomes
\begin{align}\label{10-Proof-2}
\int_M|\nabla \mathring{\rm R}_{ij}|^2=&\int_M\Big(-\frac{n}{n-2}\mathring{\rm R}_{ij}\mathring{\rm R}_{jk}
\mathring{\rm R}_{ki}-\frac{1}{n-1}R|\mathring{\rm R}_{ij}|^2\notag\\
&+\frac{(n-2)^2}{4n(n-1)}|\nabla R|^2+\frac{1}{2}|C_{ijk}|^2\Big),
\end{align}
respectively. Thus, combining \eqref{10-Proof-1} with \eqref{10-Proof-2}, we obtain
\begin{align}\label{10-Proof-3}
0=&\frac{(n-4)[(n-2)+4s]}{n-2}\int_M|\nabla \mathring{\rm R}_{ij}|^2\notag\\
&-\frac{2n[(n-1)(n-2)t+2s+(n-2)]}{(n-1)(n-2)}\int_MR|\mathring{\rm R}_{ij}|^2\notag\\
&-\frac{(n-2)[2n(n-1)t+4(n-2)s+(n^2-3n+4)]}{2n(n-1)}\int_M|\nabla R|^2\notag\\
&+\frac{2[(n-2)+(n^2-3n+4)s]}{n-2}\int_M|C_{ijk}|^2.
\end{align}

For $n\geq4$, from \eqref{2-Sec-8} we have $C_{ijk}=0$ coming from $W_{ijkl}=0$.
In particular, when $n=4$, \eqref{10-Proof-3} becomes
\begin{align}\label{10-Proof-4}
0=&(3t+s+1)\int_M(4R|\mathring{\rm R}_{ij}|^2+|\nabla R|^2),
\end{align}
which shows that if $3t+s+1\neq0$, then we have $\mathring{\rm R}_{ij}=0$ and hence $M^4$ is Einstein. This combining with
\eqref{2-Sec-4} gives that $M^4$ is of positive constant sectional curvature.

When $n\geq5$, if $t,s$ satisfy \eqref{11th-Form-1}, then \eqref{10-Proof-3} yields
\begin{align}\label{10-Proof-5}
0=&\frac{(n-4)[(n-2)+4s]}{n-2}\int_M|\nabla \mathring{\rm R}_{ij}|^2\notag\\
&-\frac{2n[(n-1)(n-2)t+2s+(n-2)]}{(n-1)(n-2)}\int_MR|\mathring{\rm R}_{ij}|^2\notag\\
&-\frac{(n-2)[2n(n-1)t+4(n-2)s+(n^2-3n+4)]}{2n(n-1)}\int_M|\nabla R|^2\notag\\
\geq&0,
\end{align}
which concludes that $M^{n}$ is Einstein. Similarly, if \eqref{11th-Form-2} is satisfied, we also have that
$M^{n}$ is Einstein and hence $M^n$ is of positive constant sectional curvature.

\subsection{Proof of Theorems \ref{9thm1} and \ref{10thm1}}
When $n=3$, \eqref{10-Proof-3} becomes
\begin{align}\label{9-Proof-3}
0=&(1+4s)\int_M|\nabla \mathring{\rm R}_{ij}|^2+3(2t+2s+1)\int_MR|\mathring{\rm R}_{ij}|^2\notag\\
&+\frac{3t+s+1}{3}\int_M|\nabla R|^2-2(1+4s)\int_M|C_{ijk}|^2.
\end{align}

If $C_{ijk,i}=0$, then \eqref{7-Proof-1} shows that $C_{ijk}=0$ and \eqref{9-Proof-3} becomes
\begin{align}\label{9-Proof-4}
0=&(1+4s)\int_M|\nabla \mathring{\rm R}_{ij}|^2+3(2t+2s+1)\int_MR|\mathring{\rm R}_{ij}|^2\notag\\
&+\frac{3t+s+1}{3}\int_M|\nabla R|^2.
\end{align}
Therefore, if $t,s$ satisfy \eqref{9th-Form-1} or \eqref{9th-Form-2}, we have $\mathring{\rm R}_{ij}=0$ and hence $M^3$ is of positive constant sectional curvature. The proof of Theorem \ref{9thm1} is finished.

If $s=-\frac{1}{4}$, then \eqref{9-Proof-3} becomes
\begin{align}\label{9-Proof-5}
0=&\frac{3}{2}(4t+1)\int_MR|\mathring{\rm R}_{ij}|^2+\frac{1}{4}(4t+1)\int_M|\nabla R|^2,
\end{align}
which is equivalent to
\begin{align}\label{9-Proof-6}
0=&6\int_MR|\mathring{\rm R}_{ij}|^2+\int_M|\nabla R|^2
\end{align}
from $t\neq-\frac{1}{4}$. Thus, we have $\mathring{\rm R}_{ij}=0$ and $M^3$ is of positive constant sectional curvature.

We complete the proof of Theorem \ref{10thm1}.


\bibliographystyle{Plain}

\end{document}